\documentclass[11pt]{amsart}
\usepackage{amsmath,amsthm,amsfonts,amscd,amssymb}

\newcounter{mycomment}

\newtheorem{prop}{Proposition}[section]
\newtheorem{lemma}[prop]{Lemma}

\newtheorem{thm}[prop]{Theorem}

\newcommand\eee{\mathbf{e}}
\newcommand\vvv{\mathbf{v}}
\newcommand\www{\mathbf{w}}
\newcommand\rrr{\mathbf{r}}
\newcommand\nnn{\mathbf{n}}

\newcommand\LLL{\mathbb{L}}
\newcommand\SSS{\mathbb{S}}
\newcommand\HHH{\mathbb{H}}

\newcommand\RRR{\mathbb{R}}

\newcommand\ddx[1]{\frac{\partial}{\partial #1}}
\usepackage{color}

\allowdisplaybreaks
\begin{document}
\title{Ruled minimal surfaces in the three dimensional Heisenberg
group}

\dedicatory{To the memory of Professor Seok Woo Kim}


\author[Y. W. Kim]{Young Wook Kim}
\address{Dept. of Mathematics, Korea University, Seoul, Korea 136-701}
\email{ywkim@korea.ac.kr}

\author[S.-E. Koh]{Sung-Eun Koh}
\address{Dept. of Mathematics, Konkuk University, Seoul, 143-701, Korea,}
\email{sekoh@konkuk.ac.kr}

\author[H. Y. Lee]{Hyung Yong Lee}
\address{Dept. of Mathematics, Korea University, Seoul, 143-701, Korea,}
\email{distgeo@korea.ac.kr}

\author[H. Shin]{Heayong Shin}
\address{Dept. of Mathematics, Chung-Ang University, Seoul, 156-756, Korea,}
\email{hshin@cau.ac.kr}

\author[S.-D. Yang]{Seong-Deog Yang}
\address{Dept. of Mathematics, Korea University, Seoul, Korea 136-701}
\email{sdyang@korea.ac.kr}

\keywords{Heisenberg group, ruled surface, minimal surface}
\subjclass[2000]{53A35 }
\date{\today}
\thanks{
The second named author was supported by KRF-2008-313-C00069.}


\begin{abstract}
It is shown that parts of planes, helicoids and hyperbolic
paraboloids are the only minimal surfaces ruled by geodesics in
the three dimensional Riemannian Heisenberg group.
 It is also shown that they are the only surfaces in the three
dimensional Heisenberg group whose mean curvature is zero with
respect to both of the standard Riemannian metric and the standard
Lorentzian metric.
\end{abstract}

\maketitle

\section{Introduction}



The three dimensional Heisenberg group $\HHH_3$ is the two-step
nilpotent Lie group standardly represented in $GL_3(\RRR)$ by
$$
\left[
\begin{array}{ccc}1&x&z+\frac12xy\\0&1&y\\0&0&1\end{array}\right].
$$
We consider in this paper two left invariant metrics on $\HHH_3,$
one is Riemannian and the other Lorentzian.
Let us denote by
$\text{Nil}^3$ the 3-dimensional Heisenberg group $\HHH_3$ endowed
with the left-invariant Riemannian metric
$$
g=dx^2+dy^2+\left(dz+\frac12(ydx-xdy)\right)^2
$$
on $\RRR^3.$ The Riemannian Heisenberg group $\text{Nil}^3$ is a
three dimensional homogeneous manifold with a 4-dimensional
isometry group; hence it is the most simple 3-manifold apart from
the space-forms.
Moreover, it is a Riemannian fibration over the Euclidean plane
$\RRR^2,$ with the projection $(x,y,z)\mapsto (x,y).$

In the first part of this paper, we give a classification of all
ruled minimal surfaces in $\text{Nil}^3.$ In order for this, we
first show in Lemma \ref{lem:rulingishorizontal} that if a ruled
surface is minimal and if a ruling geodesic is not tangent to the
fibre, then the ruled surface should be horizontally ruled. That
is, its ruling geodesics are orthogonal to the fibres. In fact, it
was one of the key observations in classifying the ruled minimal
surfaces in ${\mathbb S}^2\times{\mathbb R}$ or in ${\mathbb
H}^2\times{\mathbb R}$ in our previous paper \cite{KKSY}. It turns
out that this fact simplifies the nonlinear partial differential
equations describing ruled minimal surfaces. Then we show in
Theorem \ref{thm:ruledmin} that any ruled minimal surface in
$\text{Nil}^3$ is, up to isometries, a part of the horizontal
plane $z=0$, the vertical plane $y=0$, a helicoid $\tan(\lambda
z)=\frac{y}{x}, \lambda\ne0$ or a hyperbolic paraboloid
$z=-\frac{xy}2,$ see \S 2.3 for the definition of planes.
Moreover, we show in \S 2.8 that all of them can be regarded as
helicoids or the limits of the sequences of helicoids in the
Gromov-Hausdorff sense.

In fact, it was shown in \cite{BS} that, up to isometries,
parts of planes, the helicoids 
and the
hyperbolic paraboloids
are the only minimal
surfaces in $\text{Nil}^3$ ruled by straight lines which are
geodesics. According to Lemma \ref{lem:rulingishorizontal}, any
ruling geodesic of a ruled minimal surface is either parallel or
orthogonal to the fibres. We then note that geodesics parallel or
orthogonal to the fibres everywhere are straight lines (in the
Euclidean sense) in Lemma 2.4 and thereby show that \lq\lq
straight line" condition may be deleted in the aforementioned
claim.
For the properties of the Gauss map and representation formulae of
the minimal surfaces in $\text{Nil}^3,$ see for example
\cite{BBBI}, \cite{D}, \cite{I1}, \cite{I2}, \cite{MMP}, \cite{S}.


In the second part, we consider the natural left invariant
Lorentzian metric
$$
g_L=dx^2+dy^2-\left(dz+\frac12(ydx-xdy)\right)^2
$$
on $\HHH_3.$ (Lorentzian metrics on $\HHH_3$ are discussed in
\cite{R}, \cite{RR}). Then we consider surfaces in $\HHH_3$ whose
mean curvature is zero with respect to both metrics $g$ and $g_L$
and show that they must be one of the above mentioned surfaces,
that is, a part of planes, helicoids or hyperbolic paraboloids in
Theorem 3.2. It can be considered as a generalization of the fact
that the helicoids are the only surfaces except the planes in
$\RRR^3$ whose mean curvature is zero with respect to both the
standard Riemannian metric and the standard Lorentzian metric
\cite{K} and the fact that the helicoids (surfaces invariant under
the screw motion) are the only surfaces except the trivial ones in
$\SSS^2\times\RRR$ or $\HHH^2\times\RRR$ whose mean curvature is
zero with respect to both the standard Riemannian metric and the
standard Lorentzian metric \cite{KKSY}. In order for this we
derive the equation for the mean curvature of a graph in $\HHH_3$
to be zero with respect to the Lorentzian metric $g_L$ and compare
it with the minimal surface equation. We would like to remark that
the idea of considering these two equations in the same time is
not new, see also for example, \cite{AA}, \cite{AB}, \cite{K}.


\section{Ruled minimal surfaces in $\text{Nil}^3$}

 We first state several facts on the geometry of
$\text{Nil}^3,$ necessary for the proof of the main result in this
section. For their proofs, one may refer, for example, to
\cite{IKOS}.

\subsection{A Frame Field} It can be easily seen that
$$
\eee_1=\frac{\partial}{\partial
x}-\frac{y}2\frac{\partial}{\partial z},\quad
\eee_2=\frac{\partial}{\partial
y}+\frac{x}2\frac{\partial}{\partial z},\quad
\eee_3=\frac{\partial}{\partial z}
$$
is a left invariant orthonormal frame field on $\text{Nil}^3$ and
in particular, $\eee_3$ is tangent to the fibres. Let $\nabla$ be
the Levi-Civita connection on $\text{Nil}^3,$ then, for this frame
field we have,
\[
\begin{split}
& \nabla_{\eee_i}\eee_i=0,\quad i=1,2,3,\\
& \nabla_{\eee_1}\eee_2=
-\nabla_{\eee_2}\eee_1=\frac12\eee_3,\quad \nabla_{\eee_1}\eee_3=
\nabla_{\eee_3}\eee_1=-\frac12\eee_2,\quad \nabla_{\eee_2}\eee_3=
\nabla_{\eee_3}\eee_2=\frac12\eee_1.
\end{split}
\]

\subsection{Isometries} The isometry group of $\text{Nil}^3$ has two connected components:
an isometry either preserves the orientation of both the fibres
and the base of the fibration, or reverses both orientations. The
identity component of the isometry group of $\text{Nil}^3$ is
isomorphic to $SO(2)\ltimes\RRR^3$ whose action is given by
\[
\begin{split}
& \left( \left[ \begin{array}{cc}\cos\theta&-\sin\theta\\
\sin\theta&\cos\theta \end{array}\right], \left[
\begin{array}{c}a\\b\\c\end{array}\right]\right)\cdot\left[\begin{array}{c}x\\y\\z\end{array}\right]\\
&=\left[\begin{array}{ccc}\cos\theta&-\sin\theta&0\\\sin\theta&\cos\theta&0\\
\frac12(a\sin\theta-b\cos\theta)&\frac12(a\cos\theta+b\sin\theta)&1\end{array}\right]
\left[
\begin{array}{c}x\\y\\z\end{array}\right]+\left[\begin{array}{c}a\\b\\c\end{array}\right],
\end{split}
\]
which shows that $\text{Nil}^3$ is a homogeneous space. In fact,
one can see that, for any point $p\in \HHH_3$ and a unit tangent
vector $\vvv$ orthogonal to $\eee_3(p),$ there exists a unique
isometry $\varphi$ such that $\varphi(p)={\bf 0},
d\varphi(\vvv)=\eee_1({\bf 0})$ and
$d\varphi(\eee_3(p))=\eee_3({\bf 0}).$ Note also that the
translations along the $z$ axis (in the Euclidean sense) are
isometries belonging to the identity component.

\subsection{Euclidean Planes} A {\it Euclidean plane} or simply a
{\it plane} is a set of points $(x, y, z)\in\HHH_3$ satisfying a
linear equation $ax+by+cz+d=0.$ It is easy to see that all the
planes except the \lq\lq vertical" planes $ax+by+d=0$ are
congruent.  In fact, every nonvertical plane $ax+by+z+d=0$ is
isometric to the \lq\lq horizontal" plane $z=0$ via, for example
the isometry
$$\left[\begin{array}{c}x\\y\\z\end{array}\right]\rightarrow
\left[\begin{array}{rrr}1&0&1\\0&1&0\\-a&-b&1\end{array}\right]\left[\begin{array}{c}x\\y\\z\end{array}\right]+
\left[\begin{array}{r}-2b\\2a\\-d\end{array}\right].$$ Moreover, a
vertical plane is not congruent to a nonvertical plane since every
isometric image of a fibre is a fibre. In fact, one can check that
a vertical plane is not isometric to a nonvertical plane by
computing their curvatures.

\subsection{A parametrization of ruled surfaces}
Let $\Sigma$ be a ruled surface in $\text{Nil}^3$ and let
$p\in\Sigma$ be a point at which $T_p \Sigma$ is transversal to
the fibre. Assume, furthermore, that  the direction of the ruling
geodesic at $p$ is not perpendicular to the fibres. Then, in a
neighborhood of $p$, we can take a tangent vector field $V$ to
$\Sigma$ which is in the direction of the ruling everywhere on the
neighborhood as
$$ V =\eta(\cos\theta\eee_1 -\sin\theta\eee_2)+\eee_3 $$
for some functions $\eta$ and $\theta$ on $\Sigma$. Since $T_p
\Sigma$ is transversal to the fibre, the unit normal vector field
$\nnn$ of $\Sigma$ is not perpendicular to $\eee_3$,
$\langle\nnn,\eee_3\rangle  \neq 0$. Then
$$
 W =\sin\theta\eee_1 + \cos\theta\eee_2
   - \frac{\langle\nnn, \sin\theta\eee_1 + \cos\theta\eee_2\rangle}
   {\langle\nnn,\eee_3\rangle} \eee_3
$$
gives another tangent vector field on $\Sigma$ which is
transversal to $V$. Now we take a parametrization $X(s,t)$ of
$\Sigma$ in the neighborhood of $p$ such that $X(s,0)$ is the
integral curve of $W$ with $X(0,0)=p$ and such that $t$ parameter
curves are the ruling geodesics with $X_t(s,0)=V(X(s,0))$. Then
$X(s,t)$ is a parametrization of the ruled surface $\Sigma$ in the
neighborhood of $p$ satisfying
\begin{equation}\label{eq:parametrization}
\begin{split}
   X_s (s,0) &=\sin \alpha (s)\eee_1 + \cos\alpha (s)\eee_2 + g(s)\eee_3 , \\
   X_t (s,0) &=h(s)(\cos\alpha (s)\eee_1 -  \sin \alpha (s)\eee_2) + \eee_3,  \\
   \nabla_{X_t} X_t &=0
\end{split}
\end{equation}
for some smooth functions $h(s)$, $\alpha (s)$ and $g(s)$.

For the parametrization $X$ satisfying the condition
(\ref{eq:parametrization}), we are to compute the functions
$X_{si}$ and $X_{ti}$ defined by
\[
   \begin{aligned}
    X_s (s,t) &= X_{s1} (s,t)\eee_1 + X_{s2} (s,t)\eee_2 + X_{s3} (s,t)\eee_3, \\
    X_t (s,t) &= X_{t1} (s,t)\eee_1 + X_{t2} (s,t)\eee_2 + X_{t3}
    (s,t)\eee_3.
   \end{aligned}
\]
Now, since $t$ parameter curves are geodesics, we have
\[  \begin{split}
    \nabla_{X_t} X_t
        &= \sum_i \frac{\partial X_{ti}}{\partial t} \eee_i
           + \sum_{i,j} X_{ti}X_{tj} \nabla_{\eee_i} \eee_j  \\
        &= \bigg( \frac{\partial X_{t1}}{\partial t} + X_{t2}X_{t3} \bigg) \eee_1
           + \bigg( \frac{\partial X_{t2}}{\partial t} - X_{t1}X_{t3} \bigg) \eee_2
           + \frac{\partial X_{t3}}{\partial t} \eee_3
         = 0 \ .
\end{split}   \]
By solving the system of equations
\[  \frac{\partial X_{t1}}{\partial t} + X_{t2}X_{t3} =0,
    \frac{\partial X_{t2}}{\partial t} - X_{t1}X_{t3} =0,
    \frac{\partial X_{t3}}{\partial t} = 0 \    \]
with the initial condition
\[  X_{t1}(s,0)=h(s)\cos\alpha (s), \
    X_{t2}(s,0)= -h(s) \sin \alpha (s), \
    X_{t3}(s,0)=1   \]
we have
\[   X_{t1} (s,t) = h(s) \cos ( t-\alpha(s)),\ \
     X_{t2} (s,t) = h(s) \sin ( t-\alpha(s)),\ \
     X_{t3} (s,t) = 1.      \]
On the other hand, since the Levi-Civita connection $\nabla$ is
torsion free, one has
$$\nabla_{X_t} X_s = \nabla_{X_s} X_t.$$
Hence we have
\[  \begin{split}
     \bigg( \frac{\partial X_{s1}}{\partial t}
              + \frac12(X_{t2}X_{s3}&+X_{t3}X_{s2}) \bigg) \eee_1
     + \bigg( \frac{\partial X_{s2}}{\partial t}
              - \frac12(X_{t1}X_{s3}+X_{t3}X_{s1}) \bigg)\eee_2  \\
     &+ \bigg( \frac{\partial X_{s3}}{\partial t}
              +\frac12(X_{t1}X_{s2} - X_{t2}X_{s1}) \bigg)\eee_3   \\
  =  \bigg( \frac{\partial X_{t1}}{\partial s}
              + \frac12(X_{s2}X_{t3}&+X_{s3}X_{t2}) \bigg) \eee_1
     + \bigg( \frac{\partial X_{t2}}{\partial s}
              - \frac12(X_{s1}X_{t3}+X_{s3}X_{t1}) \bigg)\eee_2   \\
     &+ \bigg( \frac{\partial X_{t3}}{\partial s}
              +\frac12(X_{s1}X_{t2} - X_{s2}X_{t1}) \bigg)\eee_3,
\end{split}       \]
and $X_{si}$ satisfies the equations
\[  \begin{split}
         \frac{\partial X_{s1}}{\partial t} &= \frac{\partial X_{t1}}{\partial s}
             =h'(s)\cos(t-\alpha(s)) + h(s)\alpha'(s)\sin(t-\alpha(s)), \\
     \frac{\partial X_{s2}}{\partial t} &= \frac{\partial X_{t2}}{\partial s}
             =h'(s)\sin(t-\alpha(s)) - h(s)\alpha'(s)\cos(t-\alpha(s)),  \\
     \frac{\partial X_{s3}}{\partial t}
             &= \frac{\partial X_{t3}}{\partial s} + (X_{s1}X_{t2} -
             X_{s2}X_{t1})\\
              &= h(s) \sin ( t-\alpha(s))X_{s1} - h(s) \cos ( t-\alpha(s))X_{s2}
\end{split}       \]
with the initial condition
\[  X_{s1}(s,0)= \sin \alpha (s),
    X_{s2}(s,0)= \cos\alpha (s),
    X_{t3}(s,0)= g(s).  \]
By solving these equations, we get
\[  \begin{split}
    X_{s1}(s,t)&= \sin\alpha(s)+ h'(s)\sin(t-\alpha(s))+h'(s)\sin\alpha(s)  \\
     &\qquad\qquad -h(s)\alpha'(s)\cos(t-\alpha(s))+h(s)\alpha'(s)\cos\alpha(s),\\
    X_{s2}(s,t)&= \cos\alpha(s)- h'(s)\cos(t-\alpha(s))+h'(s)\cos\alpha(s) \\
     &\qquad\qquad -h(s)\alpha'(s)\sin(t-\alpha(s))-h(s)\alpha'(s)\sin\alpha(s),\\
    X_{s3}(s,t)&= g(s)-h(s)\sin t + th(s)h'(s) - h(s)h'(s)\sin t\\
     &\qquad\qquad  + h(s)^2 \alpha'(s)- h(s)^2 \alpha'(s)\cos t.
\end{split}       \]

\subsection{The second derivatives of $X$}
We are to compute the derivatives $\nabla_{X_t} X_t, \nabla_{X_s}
X_t=\nabla_{X_t} X_s$ and $ \nabla_{X_s} X_s.$ For notational
simplicity, let us set
\[  \begin{split}
    X_{t;t} &:=\nabla_{X_t} X_t = X_{tt1} \eee_1 + X_{tt2}\eee_2 + X_{tt3} \eee_3,   \\
    X_{s;t} &:=\nabla_{X_t} X_s = X_{st1} \eee_1 + X_{st2}\eee_2 + X_{st3} \eee_3,   \\
    X_{s;s} &:=\nabla_{X_s} X_s = X_{ss1} \eee_1 + X_{ss2}\eee_2 + X_{ss3} \eee_3 .  \\
\end{split}
\]
Since $t$ parameter curves are geodesics, we have $X_{t;t}=0,$
that is,
$$  X_{tt1} = X_{tt2} = X_{tt3} = 0 .$$
From
\[  \begin{split}
    X_{s;t}&= X_{t;s} \\
      &= \bigg( \frac{\partial X_{s1}}{\partial t}
              + \frac12(X_{t2}X_{s3}+X_{t3}X_{s2}) \bigg) \eee_1
        + \bigg( \frac{\partial X_{s2}}{\partial t}
              - \frac12(X_{t1}X_{s3}+X_{t3}X_{s1}) \bigg)\eee_2  \\
      &\qquad + \bigg( \frac{\partial X_{s3}}{\partial t}
              +\frac12(X_{t1}X_{s2} - X_{t2}X_{s1}) \bigg)\eee_3,   \\
     X_{s;s}
      &= \bigg( \frac{\partial X_{s1}}{\partial s}
              + X_{s2}X_{s3} \bigg) \eee_1
        + \bigg( \frac{\partial X_{s2}}{\partial s}
              - X_{s1}X_{s3} \bigg)\eee_2
        + \frac{\partial X_{s3}}{\partial s}  \eee_3
\end{split}
\]
we have 
\[ \begin{split}
   X_{st1} &= \frac{1}{2} \bigg[ \cos \alpha(s) +h'(s)\cos(t-\alpha (s))
              + h'(s) \cos\alpha (s)\\
      &\qquad\qquad\qquad\qquad\qquad\qquad\qquad
         +h(s)\alpha'(s)\sin(t-\alpha (s))-h(s)\alpha '(s) \sin \alpha (s) \\
       &\ \ +h(s) \sin(t-\alpha (s)) \bigg(g(s)+h(s)
           \bigg( -\sin t +h'(s)(t-\sin t) +2 h(s)\alpha'(s)\sin^2t/2
           \bigg) \bigg) \bigg], \\
   X_{st2} &= \frac{1}{2} \bigg[ -\sin \alpha(s)+h'(s)\sin(t-\alpha (s))
              -h'(s) \sin \alpha (s) \\
       &\qquad\qquad\qquad\qquad\qquad\qquad\qquad
          -h(s)\alpha '(s)\cos(t-\alpha (s)) -h(s) \alpha '(s)\cos\alpha (s)\\
       &\ \ + h(s)\cos(t-\alpha (s))  \bigg(-g(s)+ h(s)
         \bigg(\sin t -h'(s)(t-\sin t) -2 h(s) \alpha '(s)\sin ^2 t/2 \bigg)
         \bigg)\bigg],   \\
   X_{st3} &=\frac{1}{2} h(s) \bigg[-\cos t-h'(s)(\cos t-1)
       +h(s)\alpha '(s) \sin t \bigg],    \\
   X_{ss1} &=
      \alpha'(s)\cos\alpha(s)-2h'(s)\alpha'(s)\cos(t-\alpha(s))
                 +2h'(s)\alpha'(s)\cos\alpha(s) \\
     &\qquad\qquad\qquad\qquad\qquad\qquad\qquad
         -h(s)\alpha '(s)^2 \sin(t-\alpha(s)) -h(s)\alpha '(s)^2 \sin \alpha (s)\\
     &\ \ +\bigg(-\cos\alpha (s)+\cos (t-\alpha (s))-h'(s)\cos\alpha (s)
         +h(s)\sin (t-\alpha (s))+\alpha '(s)\sin \alpha (s) \bigg)\\
     &\qquad\qquad\qquad\qquad
           \bigg( -g(s) + h(s) \big( \sin t +h'(s)(\sin t-t)
               -2 h(s) \alpha'(s) \sin ^2t/2 \big)\bigg)\\
     &\ \ +h''(s)\sin(t-\alpha (s)) +h''(s)\sin\alpha (s)
         -h(s)\alpha''(s)\cos(t-\alpha(s))+h(s)\alpha ''(s)\cos\alpha(s) , \\
   X_{ss2} &=
      -\alpha '(s)\sin \alpha (s) -2h'(s)\alpha'(s)\sin(t-\alpha(s))
                -2h'(s)\alpha'(s)\sin\alpha(s) \\
     &\qquad\qquad\qquad\qquad\qquad\qquad\qquad
        +h(s)\alpha '(s)^2\cos (t-\alpha (s)) -h(s)\alpha '(s)^2\cos\alpha (s) \\
     &\ \ +\bigg(\sin \alpha (s)+h'(s)\big(\sin (t-\alpha(s))+\sin\alpha(s)\big)
       +2 h(s)\alpha'(s) \sin t/2 \sin (t/2-\alpha(s))\bigg)\\
     &\qquad\qquad\qquad\qquad
         \bigg( -g(s)+h(s)\big( \sin t + h'(s)(\sin t -t)
                -2 h(s) \alpha '(s) \sin^2t/2 \big)\bigg)\\
     &\ \ -h''(s)\cos (t-\alpha (s)) +h''(s)\cos\alpha (s)
        -h(s)\alpha''(s)\sin(t-\alpha(s))-h(s)\alpha''(s)\sin\alpha(s), \\
  X_{ss3} &=g'(s)+ h'(s)^2(t-\sin t) -h'(s)\left( \sin t
        -4h(s)\alpha '(s)\sin^2t/2 \right) \\
      &\qquad\qquad\qquad\qquad
         \ +h(s)\left(h''(s)(t-\sin t)-h(s)\alpha ''(s)(\cos t-1) \right).
\end{split}
\]


\subsection{Mean curvature} We give a condition for the
ruled surface $\Sigma$ to be minimal in terms of the
parametrization $X.$ Now let $E, F, G$ be the coefficients of the
first fundamental form and $l, m, n$ those of the second
fundamental form of the surface $\Sigma$ whose parametrization
satisfies (1). Then the mean curvature of $\Sigma$ in a
neighborhood of $p$ is given by
\[
\begin{split}
   H &=\frac12 \frac {Gl-2Fm+En}{EG-F^2}\\
     &=\frac12 \frac{ \langle X_t , X_t \rangle \langle X_{s;s}, X_s \times X_t
     \rangle
           -2\langle X_s , X_t \rangle\langle X_{s;t}, X_s \times X_t \rangle} {\| X_s \times
           X_t\|^3}.
\end{split}
\]
Since $$X_s \times X_t = (X_{s2}X_{t3}-X_{s3}X_{t2}) \eee_1
                   + (X_{s3}X_{t1}-X_{s1}X_{t3})\eee_2
                   + (X_{s1}X_{t2}-X_{s2}X_{t1})\eee_3 ,$$
$X$ is a parametrization of a minimal surface if and only if
\begin{equation}\label{mineq}
\begin{split}
    \tilde H :
     &= \langle X_t , X_t\rangle\langle X_{s;s}, X_s \times X_t\rangle - 2\langle X_s , X_t\rangle\langle X_{s;t}, X_s \times X_t\rangle \\
     &= \left( \sum_i X_{ti}^2 \right) \bigg( (X_{s2}X_{t3}-X_{s3}X_{t2})X_{ss1} \\
     &\qquad\qquad \qquad\qquad
       +(X_{s3}X_{t1}-X_{s1}X_{t3})X_{ss2}+(X_{s1}X_{t2}-X_{s2}X_{t1})X_{ss3}\bigg) \\
     &\quad -2 \left(\sum_i X_{si}X_{ti}\right) \bigg( (X_{s2}X_{t3}-X_{s3}X_{t2})X_{st1} \\
     &\qquad\qquad \qquad\qquad
       +(X_{s3}X_{t1}-X_{s1}X_{t3})X_{st2} +(X_{s1}X_{t2}-X_{s2}X_{t1})X_{st3}\bigg) \\
     &=0
\end{split}
\end{equation}

\subsection{Ruled minimal surfaces in $\text{Nil}^3$}

Now we are to find all ruled minimal surfaces in $\text{Nil}^3$.

\begin{lemma}\label{lem:rulingishorizontal}
If the surface whose parametrization $X$ satisfies {\rm
(\ref{eq:parametrization})} is minimal, then $h(s)=0$ for all $s$.
\end{lemma}

\noindent[Proof] Considering the parametrizations
$\tilde{X}(s,t):=X(s-s_0,t)$ if necessary, we need only to prove
$h(0)=0.$ By rotating the surface in $\text{Nil}^3$ if necessary,
we may assume that $\alpha(0)=0$. Since we have explicit formulae
for all $X_s$, $X_t$, $X_{s;s}$, $X_{s;t}$, $X_{t;t},$ we can
compute $\tilde H $ directly. In particular, since $X$ is minimal,
we have $\tilde H (0,t)=0 $ for all $t$. Since $\alpha(0)=0$,
$\tilde H (0,t) $ becomes
\[   \begin{split}
  \tilde H (0,t) = &A_0 + A_1 t + A_2 t^2 + A_3 t^3 \\
  &+B_0 \cos t +B_1 t\cos t +B_2 t^2\cos t +B_3 \cos2t +B_4 t\cos2t + B_5\cos3t\\
  &+C_0 \sin t +C_1 t\sin t +C_2 t^2\sin t +C_3 \sin2t +C_4 t\sin2t + C_5\sin3t
\end{split}     \]
where the constants $A_i$, $B_i$, $C_i$ are functions of $h(0),
h'(0), h''(0), \alpha'(0), \alpha''(0)$ and $g(0), g'(0).$
In the following computation, we are to use only the following terms:
\begin{eqnarray*}
A_3 &=& h(0)^5 h'(0)^3,\\
B_1 &=& -3 h(0)h'(0)^2  -h(0)^3 h'(0)^2 -3 h(0) h'(0)^3  -h(0)^3 h'(0)^3
                -2h(0)^3  g(0) h'(0) \alpha '(0) \\
    &  & -6 g(0) h(0)^5 h'(0) \alpha '(0) -3h(0)^3 h'(0) \alpha '(0)^2
           -9 h'(0) h(0)^5 \alpha '(0)^2 -6 h(0)^7 h'(0) \alpha '(0)^2  \\
    &  & -h(0)^4 h''(0)  -h(0)^2 h''(0), \\
B_5 &=&\frac14 \bigg(
         3h(0)^4\alpha'(0) +3h(0)^6\alpha'(0) +6h(0)^4 h'(0)\alpha'(0)
        +6h(0)^6 h'(0)\alpha'(0)    \\
    &  &\qquad\qquad   +3h(0)^4 h'(0)^2\alpha'(0)+3 h'(0)^2 \alpha '(0) h(0)^6
        -h(0)^6\alpha '(0)^3 -h(0)^8\alpha '(0)^3
    \bigg), \\
C_5 &=&\frac14\bigg(
        h(0)^3 +h(0)^5  +3h(0)^3 h'(0)  +3h(0)^5 h'(0) +3h(0)^3 h'(0)^2
       +3h(0)^5 h'(0)^2     \\
    &  &\quad  +h(0)^3 h'(0)^3  +h(0)^5 h'(0)^3
        -3h(0)^5\alpha '(0)^2  -3h(0)^7\alpha '(0)^2 -3h(0)^5 h'(0)\alpha'(0)^2    \\
    &  &\qquad\qquad \qquad\qquad\qquad\qquad\qquad\qquad\qquad\qquad
       -3 h'(0) h(0)^7\alpha '(0)^2
      \bigg).
\end{eqnarray*}
Since $\tilde H (0,t)=0 $ for all $t$ and since the above
expression is a linear combination of linearly independent
functions of $t$, all of $A_i$, $B_i$, $C_i$ must be $0$. Now from
$A_3 = h(0)^5 h'(0)^3 =0,$ we have either $h(0)=0 $ or  $h'(0)=0$.
Now suppose $h(0)\neq 0.$ Then  $h'(0)=0$ and $B_1$ becomes
\[
  B_1 = -h''(0) h(0)^4 -h''(0) h(0)^2 =-h''(0) h(0)^2(h(0)^2 + 1)= 0 .
\]
Hence we have $h''(0)= 0$ and  in addition
\[   \begin{split}
  4B_5 &=-\alpha'(0)^3 h(0)^8-\alpha'(0)^3 h(0)^6+3 \alpha'(0) h(0)^6 +3 \alpha'(0) h(0)^4 =0\\
  4C_5 &=-3 \alpha'(0)^2 h(0)^7-3 \alpha'(0)^2 h(0)^5+h(0)^5 +h(0)^3 =0  .
\end{split}     \]
Then, since
  $$ 3B_5 -  h(0) \alpha'(0) C_5 =2 \alpha'(0) h(0)^4(h(0)^2 +1)= 0 ,$$
we have $\alpha'(0) =0 $ and $C_5$ becomes
 $$ 4C_5 =h(0)^3(h(0)^2+1) =0 \ .$$
This contradicts the assumption $h(0)\ne0$. Hence we must have
$h(0)=0$ if $X$ is a parametrization of a minimal surface. \qed



\medskip

If $p$ is a point in a ruled surface $\Sigma$ at which $T_p
\Sigma$ is transversal to the fibre and the direction of the
ruling is not perpendicular to the fibres, then $\Sigma$ has the
parametrization of the type given in (\ref{eq:parametrization}) in
a neighborhood of $p$.  If, in addition, $\Sigma$ is minimal then
the above lemma implies that the direction of the ruling at $p$ is
parallel to the fibres.  This contradicts the fact that $T_p
\Sigma$ is transversal to the fibres. Therefore we can conclude
that in a ruled minimal surface $\Sigma$ the directions of the
rulings are horizontal, that is, perpendicular to the fibres
wherever $T_p \Sigma$ is transversal to the fibres.

\medskip

Now we consider the minimal surfaces which are ruled by horizontal
geodesics.

\begin{lemma}\label{lem:horzontalyruled}
If $\Sigma$ is a minimal surface in $\text{\rm Nil}^3$ ruled by
geodesics perpendicular to the  fibres, then up to the isometries
in $\text{\rm Nil}^3$, $\Sigma$ is a part of the horizontal plane
$z=0$, the vertical plane $y=0,$ a helicoid $\tan(\lambda
z)=\frac{y}{x}, \lambda\ne 0$ or a hyperbolic paraboloid
$z=-\frac{xy}2.$
\end{lemma}

\noindent[Proof] One can see that the surface $\Sigma$ has a local
parametrization $Y(s,t)$ satisfying
\begin{equation}\label{eq:horizparam}
\begin{split}
   Y_s (s,0) &=\cos\beta(s)(-\sin \alpha (s)\eee_1 + \cos\alpha (s)\eee_2)
               + \sin\beta(s)\eee_3 , \\
   Y_t (s,0) &=\cos\alpha(s)\eee_1 + \sin \alpha (s)\eee_2,   \\
   \nabla_{Y_t} Y_t &=0.
\end{split}
\end{equation}
If we set
\[
   \begin{aligned}
    Y_s(s,t) &= Y_{s1}(s,t)\eee_1 + Y_{s2}(s,t)\eee_2 + Y_{s3}(s,t)\eee_3, \\
    Y_t(s,t) &= Y_{t1}(s,t)\eee_1 + Y_{t2}(s,t)\eee_2 + Y_{t3}(s,t)\eee_3,
   \end{aligned}
\]
by solving the equation $\nabla_{Y_t} Y_t =0$ with the initial
condition
$$Y_t(s,0)=\cos\alpha(s)\eee_1 + \sin \alpha(s)\eee_2 $$
we have
\[   Y_{t1} (s,t) = \cos\alpha(s),\ \
     Y_{t2} (s,t) = \sin \alpha(s),\ \
     Y_{t3} (s,t) = 0.      \]
Moreover, from $\nabla_{Y_t} Y_s = \nabla_{Y_s} Y_t$, we can see
that $Y_{si}$ satisfies the equations
\[  \begin{split}
         \frac{\partial Y_{s1}}{\partial t} &= \frac{\partial Y_{t1}}{\partial s}
             =-\alpha'(s)\sin\alpha(s), \\
     \frac{\partial Y_{s2}}{\partial t} &= \frac{\partial Y_{t2}}{\partial s}
             =\alpha'(s)\cos\alpha(s),  \\
     \frac{\partial Y_{s3}}{\partial t}
             &= \frac{\partial Y_{t3}}{\partial s} + (Y_{s1}Y_{t2} - Y_{s2}Y_{t1})
              = \sin \alpha(s)Y_{s1} - \cos\alpha(s)Y_{s2}
\end{split}       \]
with the initial condition
\[  Y_{s1}(s,0)=-\cos\beta(s)\sin\alpha(s), \
    Y_{s2}(s,0)= \cos\beta(s)\cos\alpha(s), \
    Y_{s3}(s,0)= \sin\beta(s).  \]
By solving this system of equations, we get
\[  \begin{split}
    Y_{s1}(s,t)&= -\cos\beta(s)\sin\alpha(s)-t\alpha'(s)\sin\alpha(s),\\
    Y_{s2}(s,t)&=  \cos\beta(s)\cos\alpha(s)+t\alpha'(s)\cos\alpha(s),\\
    Y_{s3}(s,t)&= \sin\beta(s)-t\cos\beta(s)-\frac12 t^2\alpha'(s).
\end{split}       \]
By direct computations, we can see that the minimal surface
equation (\ref{mineq}) can be written as
\[
\beta'(s)
         +t\big(\alpha'(s)\beta'(s)\cos\beta(s)-\alpha''(s)\sin\beta(s)\big)
    +\frac{t^2}2 \big(\alpha'(s)\beta'(s)\sin\beta(s)+\alpha''(s)\cos\beta(s)\big)
    =0.
\]
Therefore we have $\beta'(s)=0$ and $\alpha''(s)=0,$ that is,
$\beta(s)=b$ and $\alpha(s)=as+c$ for some constants $a,b,c.$

When $a\neq0$, relocating the surface $\Sigma$ by an isometry in $\text{Nil}^3$, we may assume that
     $$\alpha(s)=as\ \ \text{and}\ \ Y(0,0)=\left(\frac{\cos b}a ,\ 0,\ 0\right).$$
Then, since
$\eee_1=\frac{\partial}{\partial x}-\frac{y}2\frac{\partial}{\partial z}$,
$\eee_2=\frac{\partial}{\partial y}+\frac{x}2\frac{\partial}{\partial z}$,
$\eee_3=\frac{\partial}{\partial z}$,
we have
\[  \begin{aligned}
   Y_s (s,0) = &-\cos b\sin(as)\eee_1 + \cos b\cos(as)\eee_2 + \sin b\eee_3  \\
     =&-\cos b\sin(as)\ddx{x}+ \cos b\cos(as)\ddx{y} \\
     & +\bigg(\sin b+\frac{y}2\cos b\sin(as)+ \frac{x}2\cos b\cos(as)\bigg)\ddx{z},  \\
   Y_t (s,t) = &\cos(as)\eee_1 + \sin(as)\eee_2\\
    =&\cos(as)\ddx{x}+\sin(as)\ddx{y}
       +\bigg(-\frac{y}2\cos(as)+\frac{x}2\sin(as)\bigg)\ddx{z}.
\end{aligned}   \]
Integrating the components of $Y_s(s,0)$ with initial data
$Y(0,0)=\left(\frac{\cos b}a,0,0 \right)$,  we have
\[
Y(s,0)=\left(\frac1{a}\cos b\cos(as) , \ \ \frac1{a}\cos
b\sin(as),
          \ \ \frac{s}{4a}(1+\cos(2b)+4a\sin b) \right).
\]
Then integrating the components of $Y_t(s,t)$ with initial data
$Y(s,0)$, we have
\[
Y(s,t)=\left(t\cos(as)+\frac1{a}\cos b\cos(as), \ \
t\sin(as)+\frac1{a}\cos b\sin(as),
          \ \ \frac{s}{4a}(1+\cos(2b)+4a\sin b)\right).
\]
Noting that
\[  Y(s,t)
  =\left(t\cos(as),\ \ t\sin(as),\ \ \frac{s}{4a}(1+\cos(2b)+4\sin b)\right),
\]
we can see that $Y$ is a parametrization of either the helicoid
$$ \tan{\lambda z}=\frac{y}{x} \ \ \text{ where }\ \
   \lambda=\frac{4a^2}{1+\cos(2b)+4a\sin b}$$
if $1+\cos(2b)+4a\sin b\neq 0$,  or the plane $z=0$ if
$1+\cos(2b)+4a\sin b = 0$.

When $a=0$ and $\cos b\neq 0$, we may assume up to isometries that
$\alpha(s)=0$ and  $Y(0,0)=(-\tan{b},0,0)$. Then
\[  \begin{aligned}
   Y_s (s,0) &=  \cos b\eee_2 + \sin b\eee_3
     =\cos b\ddx{y} +\bigg(\sin b+ \frac{x}2\cos b\bigg)\ddx{z},  \\
   Y_t (s,t) &= \eee_1 =\ddx{x} -\frac{y}2 \ddx{z},
\end{aligned}   \]
and a similar computation as above gives
\[
Y(s,t)=\left(t-\tan b ,  s\cos b ,  -\frac12 st\cos b +
\frac12s\sin b\right)
\]
which is a parametrization of the hyperbolic paraboloid $z=-\frac{xy}2$.
When $a=0$ and $\cos b= 0$,
we have $Y_s(s,0)=\eee_3$,  $Y_t(s,t)= \eee_1 =\ddx{x} -\frac{y}2 \ddx{z}$ and
$Y(s,t)$ is a parametrization of the $xz$-plane if we set $Y(0,0)=(0,0,0)$.
\qed


\begin{thm}\label{thm:ruledmin}
If $\Sigma$ is a minimal surface in $\text{\rm Nil}^3$ ruled by
geodesics, then up to the isometries in $\text{\rm Nil}^3$,
$\Sigma$ is a part of the horizontal plane $z=0$, the vertical
plane $y=0$, a helicoid $\tan(\lambda z)=\frac{y}{x}, \lambda\ne0$
or a hyperbolic paraboloid $z=-\frac{xy}2.$
\end{thm}

\noindent[Proof] If there is a point $p\in \Sigma$ at which
$T_p\Sigma$ is transversal to the fibres, then $\Sigma$ is
transversal to the fibres in a neighborhood of $p$. Therefore,
from the argument following the Lemma
{\ref{lem:rulingishorizontal}}, the ruling geodesics through any
points in the neighborhood must be horizontal. Then by the Lemma
\ref{lem:horzontalyruled} the neighborhood coincides with a part
of the helicoids, the hyperbolic paraboloid or the $xy$-plane up
to the isometries in $\text{Nil}^3$. Now since the tangent spaces
at every points of these surfaces are transversal to fibres, the
whole $\Sigma$ must be a part of one of these surfaces.

  On the other hand, if the tangent space $T_p \Sigma$ is tangent
to the fibres at every point $p\in \Sigma$, then $\eee_3$ is
tangent to $\Sigma$. Relocating $\Sigma$ by an isometry of
$\text{Nil}^3$, we may assume that $(0,0,0)\in\Sigma$ and that
$\Sigma$ is tangent to the plane $y=0$ at $(0,0,0)$. So $\Sigma$
is ruled by the fibres and has a ruled parametrization $ X(s,t) =(
x(s), y(s), t)$ satisfying $x(0)=y(0)=0$, $y'(0)=0$ and $x'(0)=1$.
The mean curvature of this parametrized surface can be easily
computed to be
\[  \frac{x''(s)y'(s)-x'(s)y''(s)}{ (x'(s)^2 +y'(s)^2 )^{3/2}}.
\]
Solving the equation $x''(s)y'(s)-x'(s)y''(s)=0$ with the above initial conditions, we have $y(s)=0$ which implies that $\Sigma$ is a part of the vertical plane $y=0$.
\qed

\medskip

By the above theorem, we know that the ruled minimal surfaces in
$\text{Nil}^3$ are congruent to the surfaces given in the theorem
which are all ruled by horizontal geodesics. In fact, the vertical
plane $y=0$ is also ruled by vertical geodesics, i.e., fibres and
this is the only doubly ruled surface among the surfaces in
Theorem {\ref{thm:ruledmin}}. Noting that isometries in
$\text{Nil}^3$ always moves fibres to fibres, we can see that the
ruled minimal surfaces in $\text{Nil}^3$ always have horizontal
ruling geodesics.  \par


\subsection{Ruled minimal surfaces as a limit of helicoids}  Consider the (generic) helicoids
$$H_{\lambda}: y - x\tan (\lambda z)=0$$ and the point
$p_{\lambda}(r_{\lambda},0,0)$ on the $x$-axis, where
$r_{\lambda}=\sqrt{2/\lambda}$. The isometry which sends $x$-axis
to itself and sends the origin to $p_{\lambda}$ is given by the
formula
\[
(x,y,z)\mapsto\left(x+r_{\lambda},y,z+\frac{r_\lambda}2y \right).
\]

If we pull back $H_{\lambda}$ via this isometry, then
$p_{\lambda}$ is moved to the origin and the equation of the
pullback of $H_{\lambda}$ becomes
\[
y-(x+r_{\lambda})\tan\left(\lambda z+\frac{r_{\lambda}\lambda}2
y\right) =0.
\]
Now we multiply this equation by $r_\lambda$ then a simple
computation shows that the equation is of the form
\[
z+\frac{xy}2+O(\sqrt{\lambda})=0
\]
and as $\lambda\to0$ this converges to the equation of the ruled
minimal surface given by
\[
z+\frac{xy}2=0.
\]
This shows that the pointed helicoids $(H_{\lambda},p_{\lambda})$
converge (in the Gromov-Hausdorff sense \cite{G}) to the
exceptional ruled minimal surface $z+xy/2=0$:
$$
(H_\lambda,p_\lambda)\to \{z+xy/2=0\} \quad\text{as }\lambda\to
0+.$$ On the other hand, one can easily check that
\begin{align*}
&(H_\lambda,0)\to \text{horizontal plane}\quad
\text{as }\lambda\to \infty,\\
&(H_\lambda,0)\to\text{vertical plane}\quad \text{as }\lambda\to
0.
\end{align*}
Therefore all the ruled minimal surface in Nil$^3$ are either the
helicoids or the limits of sequences of them.

\subsection{Straight Line Geodesics} We characterize the
geodesics which are straight lines in the Euclidean sense and give
another proof of the result in \cite{BS} mentioned in the
Introduction.

\begin{prop}\label{thm:straight}
Let $\gamma(t)=(x(t), y(t), z(t))$ be a geodesic in $\text{\rm
Nil}^3.$
\begin{itemize}
\item[(1)] If $\gamma'(0)$ is perpendicular to the fibre, then
$\gamma(t)$ is a straight line everywhere perpendicular to the
fibres. \item[(2)] If $\gamma'(0)$ is parallel to the fibre, then
$\gamma(t)$ is a straight line everywhere parallel to the fibres.
\end{itemize}
\end{prop}

\noindent[Proof] Note first that
\[
\gamma'
    = x'\frac{\partial}{\partial x} + y'\frac{\partial}{\partial y}
      + z'\frac{\partial}{\partial z} =  x'\eee_1 + y'\eee_2
      + \left(z'+\frac12(x'y-xy')\right)\eee_3.
\]
Then we have
\[
\begin{aligned}
  \nabla_{\gamma'}\gamma'
      &= x'' \eee_1 + y''\eee_2
         +  \left(z'+\frac12(x'y-xy')\right)'\eee_3\\
      &\quad +x'\nabla_{\gamma'}\eee_1 +y'\nabla_{\gamma'}\eee_2
         +\left(z'+\frac12(x'y-xy')\right)\nabla_{\gamma'}\eee_3    \\
      &= \left(x''+y'(z'+\frac12(x'y-xy'))\right)\eee_1
          + \left(y''-x'(z'+\frac12(x'y-xy'))\right)\eee_2 \\
      &\quad  + \left(z'+\frac12(x'y-xy')\right)'\eee_3.
\end{aligned}
\]
Hence $\gamma(t)=(x(t), y(t), z(t))$ is a geodesic if and only if
\begin{equation}\label{eq:geodesiceq}
\begin{aligned}
&x''+y'(z'+\frac12(x'y-xy'))=0,\\
&y''-x'(z'+\frac12(x'y-xy'))=0,\\
&\left(z'+\frac12(x'y-xy')\right)'=0.
\end{aligned}
\end{equation}
Note that the straight line $(a, b, ct+d)$ parallel to the fibre
is a geodesic. Now, suppose $\langle \gamma'(0), \eee_3\rangle=0.$
Then, since
$$\langle \gamma'(0),
\eee_3\rangle=\left(z'+\frac12(x'y-xy')\right)(0)=0$$ and since
$z'+\frac12(x'y-xy')$ is a constant function from the geodesic
equation (\ref{eq:geodesiceq}), we have $z'+\frac12(x'y-xy')=0$
for all $t.$ Moreover, the geodesic equation (\ref{eq:geodesiceq})
gives
$$x''(t)=y''(t)=0,$$
that is, $x(t)$ and $y(t)$ is a linear function of $t$ and
consequently from the geodesic equation (\ref{eq:geodesiceq})
again, we have
$$z(t)= -\frac12(x'(0)y(0)-x(0)y'(0))t+c$$
for a constant $c.$ Now it is easy to see that $\gamma(t)$ is
perpendicular to the fibres everywhere.

If $\gamma'(0)$ is parallel to the fibre, then the fibre through
$\gamma(0)$ is an image of a geodesic, from the uniqueness of the
geodesic, we have $\gamma(t)=(x(0), y(0), at+b)$ for constants $a,
b$ which is parallel to the fibre everywhere. \qed

\begin{prop}
Suppose the straight line $\delta(t)=(a_1t+b_1, a_2t+b_2,
a_3t+b_3)$ is a geodesic in $\text{\rm Nil}^3.$ Then
$\delta'(0)=(a_1, a_2, a_3)$ is either perpendicular or parallel
to the fibre. Moreover, if $\delta'(0)$ is perpendicular to the
fibre, then $\delta(t)$ is perpendicular to the fibre everywhere
and if $\delta'(0)$ is parallel to the fibre, then $\delta(t)$ is
parallel to the fibre everywhere.
\end{prop}

\noindent[Proof] In the proof of the above Proposition
\ref{thm:straight}, one can see that in order for the straight
line $\delta(t)$ to be a geodesic, it should be that
$$a_3=-\frac12(a_1b_2-a_2b_1).$$ The claims follow easily form this fact. \qed

\medskip

Now we can also say that every ruled minimal surfaces in
$\text{Nil}^3$ is ruled by geodesics which are also straight
lines. We remark that it was shown in \cite{BS} that if the
surface is ruled by geodesics which are also straight lines then
the surface must be a part of the planes, helicoids or hyperbolic
paraboloids, however, in view of Theorem {\ref{thm:ruledmin}}, we
can see that the \lq\lq straight line" condition is redundant. On
the other hand, one may get Theorem {\ref{thm:ruledmin}} by
applying the aforementioned result together with Lemma
\ref{lem:rulingishorizontal} and Proposition \ref{thm:straight}.


\section{Another characterization of ruled minimal surfaces in $\HHH_3.$}

We consider surfaces in $\HHH_3$ whose mean curvature is zero with
respect to both metrics $g$ and $g_L$ and show that they must be
one of (a part of) the above mentioned surfaces, that is, planes,
helicoids and hyperbolic paraboloids.

\subsection{A Lorentzian connection}
Let us consider the left-invariant Lorentzian metric
$$ g_L=dx^2+dy^2-\left(dz+\frac12(ydx-xdy)\right)^2$$
 on $\HHH_3$
and let $\langle\cdot,\cdot\rangle$ be the Lorentzian inner
product. Let $\eee_1, \eee_2$ and $\eee_3$ be the same as the ones
given in \S 2. It is easy to show that they are orthonormal with
respect to the Lorentzian metric $g_L$ as well, that is, $\langle
\eee_i, \eee_j\rangle=0$ if $i\ne j$ and
$$
\langle\eee_1, \eee_1\rangle=\langle\eee_2, \eee_2\rangle=1,\quad
\langle\eee_3, \eee_3\rangle=-1.$$
 Now
let $D$ be the Levi-Civita connection for the metric $g_L.$

\begin{prop} We have
\[
\begin{split}
&D_{\eee_1}\eee_2=-D_{\eee_2}\eee_1=\frac12\eee_3,\quad
D_{\eee_1}\eee_3=D_{\eee_3}\eee_1=\frac12\eee_2,\\
&D_{\eee_2}\eee_3=D_{\eee_3}\eee_2=-\frac12\eee_1,\quad
D_{\eee_i}\eee_i=0, \quad i=1,2,3.
\end{split}
\]
\end{prop}

\noindent[Proof] It is known that the Koszul formula
\[
\begin{split}
2\langle \nabla_V W,X\rangle
={} &V\langle W,X\rangle + W \langle X,V\rangle - X\langle V,W\rangle\\
&-\langle V,[W,X]\rangle+\langle W,[X,V]\rangle+\langle
X,[V,W]\rangle
\end{split}
\]
holds, see, for instance, \cite{O}. Since
\[
[\eee_1,\eee_2]=\eee_3,\quad [\eee_2,\eee_3]=[\eee_3,\eee_1]=0,
\]
one has
\[
\begin{split}
&\langle D_{\eee_1}\eee_2,\eee_1\rangle=0,\\
&\langle D_{\eee_1}\eee_2,\eee_2\rangle=0,\\
&2\langle D_{\eee_1}\eee_2,\eee_3\rangle
=\langle\eee_3,[\eee_1,\eee_2]\rangle
=\langle\eee_3,\eee_3\rangle=-1
\end{split}
\]
and
\[
D_{\eee_1}\eee_2=\frac12 \eee_3.
\]
Since
\[
\begin{split}
&\langle D_{\eee_1}\eee_3,\eee_1\rangle=0,\\
&\langle D_{\eee_1}\eee_3,\eee_3\rangle=0,\\
&2\langle D_{\eee_1}\eee_3,\eee_2\rangle
=\langle\eee_3,[\eee_2,\eee_1]\rangle
=\langle\eee_3,-\eee_3\rangle=1,
\end{split}
\]
one has
\[
D_{\eee_1}\eee_3=\frac12 \eee_2.
\]
One can check the others in the same manner.\qed

\subsection{Lorentzian Exterior Product}
For tangent vectors $$\vvv=a_1\eee_1+a_2\eee_2+a_3\eee_3,
\www=b_1\eee_1+b_2\eee_2+b_3\eee_3$$ in $\text{Nil}_1^3,$ the
Lorentzian exterior product $\vvv\times_L\www$ is computed as
\begin{eqnarray*}
\vvv\times_L\www&=&\left|
\begin{array}{rrr}\eee_1&\eee_2&-\eee_3\\a_1&a_2&a_3\\b_1&b_2&b_3\end{array}\right|\\
&=&(a_2b_3-a_3b_2)\eee_1+(a_3b_1-a_1b_3)\eee_2+(a_2b_1-a_1b_2)\eee_3
\end{eqnarray*}
which is orthogonal to both $\vvv$ and $\www.$ One can easily see
that $\vvv\times_L\www=\bf{0}$ if and only if $\vvv$ and $\www$
are linearly dependent.

\subsection{Zero mean curvature equation}
Let $\Sigma$ be a graph of a function
 $z=f(x,y)$ in $\HHH_3$ and consider
 the
parametrization $\rrr(x,y)=(x,y,f(x,y))$ of $\Sigma.$ Set
\[
p=f_x+\frac{y}2,\ q=f_y-\frac{x}2.
\]
If $\Sigma$ is minimal, that is, the mean curvature is zero in
$\text{Nil}^3,$
 the function $f$ satisfies the minimal surface equation
 $$ (1+q^2)f_{xx}-2pqf_{xy}+(1+p^2)f_{yy}=0.$$
For the derivation of this equation, see for example \cite{IKOS}.

In this section, we are to derive an equation for the mean
curvature of the graph $\Sigma$ to be zero with respect to the
Lorentzian metric $g_L.$ First, let us recall some definitions. A
point $z\in \Sigma$ is called {\it spacelike} if the induced
metric on $T_z\Sigma$ is Riemannian, {\it timelike} if the induced
metric is Lorentzian and {\it lightlike} if the induced metric has
rank $1.$ We are to derive the equation when $\Sigma$ is
spacelike, that is, every point of $\Sigma$ is a spacelike point.
The case when $\Sigma$ is timelike is almost identical. Note that
when $z\in \Sigma$ is lightlike, one cannot define the mean
curvature.

Now let $\Sigma$ be a spacelike graph of a function
 $z=f(x,y).$
Note first that $p^2+q^2<1$ since the graph is spacelike. We now
compute the first fundamental form ${\rm I}$ and the second
fundamental form ${\rm II}$ of $\Sigma.$ Since
\[
\begin{split}
&\rrr_x=(1,0,f_x)=\eee_1+p\eee_3, \quad
\rrr_y=(0,1,f_y)=\eee_2+q\eee_3,
\\
&\langle\rrr_x,\rrr_x\rangle=1-p^2,\quad
\langle\rrr_x,\rrr_y\rangle=-pq,\quad
\langle\rrr_y,\rrr_y\rangle=1-q^2.
\end{split}
\]
one has
$$
E=\langle\rrr_x,\rrr_x\rangle=1-p^2,\quad
F=\langle\rrr_x,\rrr_y\rangle=-pq.\quad
G=\langle\rrr_y,\rrr_y\rangle=1-q^2. $$
Since
\[
\rrr_x\times_L\rrr_y =-p\eee_1-q\eee_2-\eee_3
\]
the unit normal vector field $\nnn$ to the graph is
\[
\nnn=\frac1{W} \left(-p\eee_1-q\eee_2-\eee_3\right), \quad
W=\sqrt{1-(p^2+q^2)}.
\]
Since the directional derivatives of $p$ and $q,$ $\eee_i(p),
\eee_i(q)$ are computed as
\[
\begin{split}
\eee_1(p) & = \left(\frac{\partial}{\partial x}-\frac{y}2\frac{\partial}{\partial z}\right)\left(f_x+\frac{y}2\right)=f_{xx},\\
\eee_1(q) & = \left(\frac{\partial}{\partial
x}-\frac{y}2\frac{\partial}{\partial z}\right)
\left(f_y-\frac{x}2\right)=f_{xy}-\frac12,\\
\eee_2(p) & = \left(\frac{\partial}{\partial
y}+\frac{x}2\frac{\partial}{\partial
z}\right)\left(f_x+\frac{y}2\right)=f_{xy}+\frac12,\\
 \eee_2(q) &
= \left(\frac{\partial}{\partial
y}+\frac{x}2\frac{\partial}{\partial z}\right)
\left(f_y-\frac{x}2\right)=f_{yy},
\end{split}
\]
one has
\begin{eqnarray*}
D_{\rrr_x}\rrr_x&=&D_{(\eee_1+p\eee_3)}(\eee_1+p\eee_3)\\
&=&p\eee_2+f_{xx}\eee_3,\\
D_{\rrr_y}\rrr_x&=& -\frac{p}2\eee_1+\frac{q}2\eee_2+f_{xy}\eee_3, \\
D_{\rrr_y}\rrr_y&=& -q\eee_1+f_{yy}\eee_3.
\end{eqnarray*}
Then one has the following coefficients of the second fundamental
form ${\rm II}.$
\begin{eqnarray*}
l &=&\langle D_{\rrr_x}\rrr_x, \nnn\rangle =
\frac1{W}\left(-pq+f_{xx}\right),\\
m &=&\langle D_{\rrr_y}\rrr_x, \nnn\rangle =
\frac1{W}\left(\frac{p^2}2-\frac{q^2}2+f_{xy}\right),\\
n &=&\langle D_{\rrr_y}\rrr_y, \nnn\rangle =
\frac1{W}\left(pq+f_{yy}\right).
\end{eqnarray*}
Now the mean curvature $H$ of the spacelike graph $\Sigma$ is
computed as
\[
H=\frac12\frac{lG-2mF+nE}{EG-F^2}.
\]
Then, since
\[
\begin{split}
lG-2mF+nE &=\frac1{W}\left[(-pq+f_{xx})(1-q^2) +
\left(\frac{p^2-q^2}2+f_{xy}\right)pq
+(pq+f_{yy})(1-p^2)\right]\\
&=\frac1{W}\left[(1-q^2)f_{xx}+2pqf_{xy}+(1-p^2)f_{yy}\right],
\end{split}
\]
one can see that the mean curvature of the graph $z=f(x,y)$ of a
function $f(x,y)$ is zero if and only if
$$(1-q^2)f_{xx}+2pqf_{xy}+(1-p^2)f_{yy}=0.$$
When the graph $\Sigma$ is timelike, one has the same equation.


%

\subsection{Zero mean curvature surface}

We prove the following theorem:

\begin{thm}\label{thm:minmax} Let $\Sigma$ be a surface in $\HHH_3.$
If the mean curvature of $\Sigma$ is zero with respect to both
metrics $g$ and $g_L,$ then up to the isometries in $\text{\rm
Nil}^3$, $\Sigma$ is a part of the horizontal plane $z=0$, the
vertical plane $y=0,$ a helicoid $\tan(\lambda z)=\frac{y}{x},
\lambda\ne0$ or a hyperbolic paraboloid $z=-\frac{xy}2.$
\end{thm}

\noindent[Proof] Suppose first that $\Sigma$ has a point around
which can be represented as a graph of a function of $(x,y),$ say,
$z=f(x,y).$ Consider the vector field
$$X=-q\eee_1+p\eee_2.$$
Since
\[ X=-q\eee_1+p\eee_2
=-q\rrr_x+p\rrr_y
\]
it is tangent to $\Sigma.$ Since the vector \[
N=\rrr_x\times\rrr_y =-p\eee_1-q\eee_2-\eee_3
\]
is orthogonal to $\Sigma$ and since $
N\times\eee_3=-q\eee_1+p\eee_2=X,$ $X$ is orthogonal to both $N$
and $\eee_3.$ Then one has
\[
\nabla_XX
= \left(q\left(f_{xy}-\frac12\right)-pf_{yy}\right)\eee_1 +
\left(p\left(f_{xy}+\frac12\right)-qf_{xx}\right)\eee_2.
\]
Now, since the mean curvature of $\Sigma\subset \HHH_3$ is zero
with respect to both $g$ and $g_L,$ one has
\begin{align}
&(1+q^2)f_{xx}-2pqf_{xy}+(1+p^2)f_{yy}=0,\label{eq:mineq}\\
&(1-q^2)f_{xx}+2pqf_{xy}+(1-p^2)f_{yy}=0.\label{eq:maxeq}
\end{align}
Subtracting two equations, one has
\begin{equation}\label{eq:minmax}
q^2f_{xx}-2pqf_{xy}+p^2f_{yy}=0
\end{equation}
and then one has finally by \eqref{eq:minmax}
\[
\begin{split}
X\times\nabla_XX&=(-q\eee_1+p\eee_2)\times\left[
\left(q\left(f_{xy}-\frac12\right)-pf_{yy}\right)\eee_1 +
\left(p\left(f_{xy}+\frac12\right)-qf_{xx}\right)\eee_2\right]\\
&=(q^2f_{xx}-2pqf_{xy}+p^2f_{yy})\eee_3\\&=0.
\end{split}
\]
Now, since $X$ and $\nabla_XX$ are of the same direction, the
integral curve of $X$ passing through a point in $\Sigma$ is a
geodesic and since $X$ is orthogonal to $\eee_3,$ the geodesic is
orthogonal to the fibre. Hence the surface $\Sigma$ is a
horizontally ruled minimal surface in $\text{Nil}^3.$

If the surface $\Sigma$ has no point around which $\Sigma$ is
represented as the graph of $f(x,y),$ then it is a vertical
cylinder over a curve in the $xy$ plane and has a parametrization
$X(s,t)=(x(s),y(s),t),\  x(0)=y(0)=0.$ By repeating the arguments
in Theorem \ref{thm:ruledmin}, one can show that the surface is
isometric to the vertical plane $y=0.$ Now this completes the
proof.
 \qed


\subsection{Remark} If we add (\ref{eq:mineq}) and (\ref{eq:maxeq}), we have
\[ f_{xx}+f_{yy}=0 \]
that is, if a graph of a function $z=f(x,y)$ in $\HHH_3$ satisfies
the condition of Theorem \ref{thm:minmax}, $f$ must be a harmonic
function. This fact is true for the three dimensional Lorentzian
space $\LLL^3$ and is the motivation of \cite{KLY}. We think it is
a nontrivial fact and would like to find applications of this fact
in the future study.

\end{document}